\documentclass{amsart}
\usepackage{tikz}
\usetikzlibrary{intersections,arrows.meta}
\usepackage{enumerate}
\usepackage{hyperref}
\usepackage{amssymb}

\newcounter{minutes}
\divide\time by 60
\newcounter{hours}
\multiply\time by 60 \addtocounter{minutes}{-\time}

\newtheorem{theorem}{Theorem}[section]
\newtheorem{lemma}[theorem]{Lemma}

\theoremstyle{definition}

\theoremstyle{remark}
\newtheorem{remark}[theorem]{Remark}

\numberwithin{equation}{section}

\newcommand{\degree}{-{1}/{2}+ \mathrm{i} \tau}

\title[Monotonicity properties of Ferrers and conical functions]{Monotonicity properties of ratios of Ferrers and conical functions related to a Riemannian nonlocal eigenvalue problem}

\author[\'A. Baricz]{\'Arp\'ad Baricz}
\address{Department of Economics, Babe\c{s}-Bolyai University, Cluj-Napoca 400591, Romania and Institute of Applied Mathematics, \'Obuda University, 1034 Budapest, Hungary}
\email{bariczocsi@yahoo.com}

\author[M.M. Rudrappa]{Manjunatha M. Rudrappa}
\address{Department of Mathematics, Indian Institute of Technology Indore, Indore 453552, India}
\email{manjunatha.m.r@icloud.com}

\author[S. Singh]{Sanjeev Singh}
\address{Department of Mathematics, Indian Institute of Technology Indore, Indore 453552, India}
\email{snjvsngh@iiti.ac.in}

\subjclass[2010]{Primary: 33C05, Secondary: 26D15, 35P15.}
\keywords{Hypergeometric functions; Ferrers functions of the first kind; conical functions of the first kind; nonlocal eigenvalue problem}

\begin{document}

\begin{abstract}
Recently, for a given Riemannian manifold with certain curvature constraints, S. Kaj\'ant\'o and A. Krist\'aly \cite{SKAK24} formulated a variational problem (on either positively or negatively curved Riemannian model spaces) and computed the eigenvalues and eigenfunctions in terms of the radii, curvature and weight. Moreover, they conjectured that a saturation phenomenon exists also in the case of positively or negatively curved Riemannian model spaces (not only in the Euclidean case), and this is based entirely on the incompatibility of a system of nonlinear equations involving special functions, like Ferrers functions of the first kind and conical functions of the first kind. In this paper, we provide a partial solution for both positive and negative curvature case.
\end{abstract}
	
\maketitle

\section{Introduction}\label{Introduction}
In this paper, our aim is to verify partially the sufficient condition of a conjecture related to a Riemannian nonlocal eigenvalue problem. In 2024, S. Kaj\'ant\'o and A. Krist\'aly \cite{SKAK24} investigated the shape optimisation problem in a complete $n$-dimensional Riemannian manifold $M_k^n$ with constant sectional curvature $k$. The investigated problem is to find the optimal domain $\Omega \subset M_k^n$ for which we obtain the minimum eigenvalue of a nonlocal eigenvalue problem
\begin{equation}\label{eq:EigenValueProblem}
	\begin{cases}
		-\Delta_k u + \alpha \displaystyle\int_{\Omega} u \, {\rm d}v_{k} = \beta u \quad \text{in } \Omega,\\[1.5ex]
		u = 0 \quad \text{on } \partial\Omega,
	\end{cases}
\end{equation}
among the domains $\Omega$ which are either a geodesic ball or a union of two disjoint geodesic balls and have measure $\omega$. Here $\Delta_k$ is the Laplace-Beltrami operator with curvature $k$. They formulated the conjecture and derived a sufficient condition for it (see \cite[Conjecture 5.2]{SKAK24} or \cite[Condition 5.5]{SKAK24}). The conjecture states that for every positive $\alpha$ less than the threshold $\alpha_{\omega,k,n}^{\mathrm{crit}}$, the optimal domain is a single geodesic ball. The problem \eqref{eq:EigenValueProblem} is a natural extension to the shape optimisation problem studied in \cite{euclidSaturation}, where $M_{k=0}^n=\mathbb{R}^{n}$, and in this case the corresponding optimisation problem is a well-established result. In order to have a partial verification of \cite[Condition 5.5]{SKAK24}, we are going to investigate the monotonicity of the ratio of the Ferrers functions of the first kind $y_1(x)={\mathsf P}^{-\mu}_{\nu}(x)$ and $y_2(x)={\mathsf P}^{1-\mu}_{\nu}(x)$ on the interval $(-1, 1)$ and the ratio of the conical (or Mehler) functions of the first kind $y_3(x)=P^{-\mu}_{\degree}(x)$ and $y_4(x)=P^{1-\mu}_{\degree}(x)$ on the interval $(1, \infty)$. It is worth mentioning here that the function $x \mapsto {\mathsf P}^{-\mu}_{\nu}(x)$ is a real-valued solution with real order $-\mu$ and degree $\nu$ of the associated Legendre differential equation
\begin{equation}\label{eq:legendreequation}
	(1 - x^2)y''(x) - 2xy'(x) + \left[ \nu (\nu + 1) - \frac{\mu^2}{1 - x^2} \right]y(x) = 0,
\end{equation}
for $x \in \left(-1,1\right)$. Similarly, the function $x \mapsto P^{-\mu}_{\degree}(x)$ is a real-valued solution with real order $-\mu$ and complex degree $\nu=\degree$ of the associated Legendre equation \eqref{eq:legendreequation} on $\left(1,\infty\right)$ where $\tau \in \mathbb{R}$. Recall that the number of zeros of $y_1$ and $y_2$ is finite in $\left(-1,1\right)$ for every real $\mu$ and $\nu$. Throughout this article, let $N_1$ and $N_2$ denote the number of zeros of $y_1$ and $y_2$, respectively. Let $x_{1, 1}<x_{1,2}<\ldots<x_{1,N_1}$ be the zeros of $y_1$ and $x_{2, 1}<\ldots<x_{2, N_2}$ be the zeros of $y_2$. If $\tau=0,$ then both functions $y_3$ and $y_4$ become associated Legendre functions of the first kind, defined on $(1,\infty),$ see \cite[p. 353]{Nist}. However, this case is irrelevant in the study of Condition 5.5, mentioned above. Hence, we assume $\tau \neq 0$. The number of zeros of $y_3$ and $y_4$ is infinite whenever $\tau \neq 0$. This is a direct consequence of Leighton's lemma \cite[Theorem 31.7]{AgarwalORegan} with a suitable transformation. For completeness, we have provided this in Lemma \ref{lemma:infintezeros}. Throughout this paper, let $x_{3, 1}<x_{3,2}<\dots$ and $x_{4, 1}<x_{4,2}<\ldots$ be the zeros of $y_3$ and $y_4$.

The paper is organized as follows: the next section contains the monotonicity of the ratios $y_1/y_2$ and $y_3/y_4,$ together with the partial verification of Condition 5.5, while section 3 contains the proofs of the main results contained in section 2.
	
\section{Preliminaries and main results}
\label{PrelANDResulst}

We first write the definitions of the volume function $V$, the profile function $f$, and the ratio function $l$ from \cite{SKAK24}. The values of  $r \mapsto V(r)$, $r \mapsto f(r)$ and $r \mapsto l(r)$ are given as follows
\begin{equation}
	V(r) = \begin{cases}
		n\omega_n \displaystyle\int_0^r \left(\frac{\sin(\sqrt{k}t)}{\sqrt{k}}\right)^{n-1} {\rm d}t, \text{ for } r \in \left[0, \frac{\pi}{\sqrt{k}}\right],\ k > 0 \\
		n\omega_n \displaystyle\int_0^r \left(\frac{\sinh(\sqrt{-k}t)}{\sqrt{-k}}\right)^{n-1} {\rm d}t, \text{ for } r \in [0, \infty),\ k < 0 \end{cases},
\end{equation}
\begin{equation}
	f(r) = \begin{cases}
		\left(\sin(\sqrt{k}r)\right)^{1-\frac{n}{2}} {\mathsf P}^{1-\frac{n}{2}}_{\gamma_{\beta,k,n}}\left(\cos(\sqrt{k}r)\right),  \text{ for } r \in \left(0, \frac{\pi}{\sqrt{k}}\right),\ k > 0 \\[15pt]
		\left(\sinh(\sqrt{-k}r)\right)^{1-\frac{n}{2}}P^{1-\frac{n}{2}}_{\gamma_{\beta,k,n}}\left(\cosh(\sqrt{-k}r)\right),  \text{ for } r \in (0, \infty),\ k < 0
	\end{cases}
\end{equation}
and
\begin{equation}
	l(r) = \begin{cases}
		\left(\sin(\sqrt{k}r)\right)^{n-1} \left(\displaystyle\frac{{\mathsf P}^{-\frac{n}{2}}_{\gamma_{\beta,k,n}}(\cos(\sqrt{k}r))}{{\mathsf P}^{1-\frac{n}{2}}_{\gamma_{\beta,k,n}}(\cos(\sqrt{k}r))}\right)^{2},  \text{ for } r \in \left(0, \frac{\pi}{\sqrt{k}}\right),\ k > 0 \\[15pt]
	\left(\sinh(\sqrt{-k}r)\right)^{n-1}\left(\displaystyle\frac{P^{-\frac{n}{2}}_{\gamma_{\beta,k,n}}(\cosh(\sqrt{-k}r))}{P^{1-\frac{n}{2}}_{\gamma_{\beta,k,n}}(\cosh(\sqrt{-k}r))}\right)^{2},  \text{ for } r \in (0, \infty),\ k < 0
	\end{cases}.
\end{equation}
Here $\omega_n$ is the volume of the unit $n$-dimensional Euclidean ball, and the degree of Ferrers or conical functions is
$$\gamma_{\beta,k,n} = -\frac{1}{2}+\sqrt{\frac{(n-1)^2}{4}+\frac{\beta}{k}},$$ where $\beta>0$ is the first eigenvalue of the local eigenvalue problem \cite[eq. 12]{SKAK24} on the geodesic ball $B(\tilde{R})$ of radius $\tilde{R}$. Now we restate Condition 5.5.

{\bf Condition 5.5}. For every $R \in \left(\tilde{R}, \tilde{R}_2\right)$, the system of equations
$$l(r_1)=l(r_2) \quad \text{and}\quad V(r_1)+V(r_2)=V(R)$$
is incompatible. Here $\tilde{R}$ is the first zero of a profile function $f$ and $\tilde{R}_2$ is such a number that $V(\tilde{R}_2)= 2 V(\tilde{R})$ and $0 \leq r_1<\tilde{R}<r_2 \leq R$.

If $k>0$, then $\gamma_{\beta,k,n}$ is real, and we have to deal with the Ferrers function of the first kind. If $k<0$, then $\gamma_{\beta,k,n}$ is of the form $-1/2+ \mathrm{i} \tau_{\beta,k,n}$, where $\tau_{\beta,k,n}>0$ and is defined as
\begin{equation}\label{eqtaudef}\tau_{\beta,k,n}=\sqrt{\frac{-(n-1)^2}{4}-\frac{\beta}{k}}\end{equation}
and we deal with conical functions of the first kind. Now, we can write
\begin{equation}
	l(r) = \begin{cases}
		(1-x^2)^{{\frac{n-1}{2}} } \left( \frac{y_1(x)}{y_2(x)} \right)^2, \text{ for } x=\cos(\sqrt{k}r)\in(-1,1),\ k>0 \\
		(x^2-1)^{{\frac{n-1}{2}} } \left( \frac{y_3(x)}{y_4(x)} \right)^2, \text{ for } x=\cosh(\sqrt{-k}r)\in (1, \infty),\ k<0	
	\end{cases}.
\end{equation}
Since the function $r \mapsto x=\cos(\sqrt{k}r)$ is decreasing on $\left[0, \pi/\sqrt{k}\right]$, the first zero $\tilde{R}$ of $f$ corresponds to the last zero in the $x$ domain whenever $k>0$. However, $r \mapsto x=\cosh(\sqrt{-k}r)$ is increasing on $(0, \infty)$. Hence, if $k<0$, then the first zero $\tilde{R}$ of $f$ corresponds to the first zero in the $x$ domain. Our key tool towards verification of the above mentioned Condition 5.5 is to use the monotonicity of $l$ presented in Theorems \ref{thm:l-monotonicityk>0} and \ref{thm:l-monotonicityk<0}. To achieve this result, first we obtain the monotonicity of $x \mapsto y_1(x)/y_2(x)$ and $x \mapsto y_3(x)/y_4(x)$ between two consecutive zeros of the denominators for any $\mu, \nu \in \mathbb{R}$ and $\tau >0$.

To obtain these monotonicity properties, we use the interlacing property of the zeros of the corresponding function pairs, namely $(y_1, y_2)$ and $(y_3, y_4)$. This interlacing behaviour is a direct consequence of the functions in each pair having only simple zeros, sharing no zeros and satisfying the system of difference-differential equations pair, that is
\begin{equation} \label{eq:y1coupledpart}
	(1-x^2)y_{1}'(x) = -(1-x^2)^{1/2}y_{2}(x) + \mu x y_{1}(x),
\end{equation}
\begin{equation}\label{eq:y2coupledpart}
	(1-x^2)y_{2}'(x) = (1-\mu) x y_{2}(x) + (\mu+\nu)( \nu-\mu+1) (1-x^2)^{1/2} y_{1}(x)
\end{equation}
and
\begin{equation} \label{eq:y3coupledpart}
	(1-x^2)y_3'(x) = -(x^2-1)^{1/2}y_{4}(x) + \mu x y_{3}(x),
\end{equation}
\begin{equation}\label{eq:y4coupledpart}
	(1-x^2)y_4'(x) = (1-\mu) x y_{4}(x) +((\mu -1/2)^2 +\tau ^2) (x^2 -1 )^{1/2} y_{3}(x).
\end{equation}
For equations \eqref{eq:y1coupledpart} and \eqref{eq:y2coupledpart} we refer to \cite[eq. 14.10.1, eq. 14.10.2  and eq .14.10.4]{Nist}. Similarly \eqref{eq:y3coupledpart} and \eqref{eq:y4coupledpart} are due to \cite[eq. 14.10.4, eq. 14.10.6  and eq .14.10.7]{Nist}. For more details about conical functions, see  \cite[p. 375]{Nist}. For the sake of completeness, we have given these properties of zeros in Lemmas \ref{lemma:simplezeros}, \ref{lemma:NoCommonZeros} and \ref{lemma:interlace}. We now first study the monotonicity of the ratio of Ferrers functions and $l$ to proceed to Condition 5.5 in the case of positive curvature. Then we do the same for the ratio of conical functions, which is associated with negative curvature.

\subsection{Monotocity of the ratio $y_1/y_2$ involving Ferrers functions}

First observe that in the verification of Condition 5.5 for $k>0$ we have $\nu =\gamma_{\beta,k,n} \in \mathbb{R}$, $-\mu=-n/2<0$ and $1-\mu \leq 0$. Thus in view of the distribution of the zeros of Ferrers functions \cite[p. 368]{Nist}, we have $$N_1=\max\left\{\left\lceil-\frac{1}{2}+\sqrt{\frac{(n-1)^2}{4}+\frac{\beta}{k}} - \frac{n}{2} \right\rceil,0\right\}$$ and $N_2= N_1+1$. Here $\lceil \epsilon \rceil$ is the integer greater than or equal to $\epsilon$. Notice that if $\epsilon>0$ then $\lceil \epsilon \rceil \geq 1$. Since $n \geq 2$, $\beta>0$ and $k>0,$ it is not difficult to see that if $\nu= \gamma_{\beta,k,n},$ then $N_1 \geq 1$ and $N_2 \geq 2$ if and only if $n < \beta/k$. Hence, in the case of Condition 5.5 with positive curvature, we have to deal with the ratio $y_1(x)/y_2(x)$ such that the denominator has at least one zero. This motivates us to study the monotonicity of the ratio of Ferrers functions, where $\mu$ and $\nu$ are real numbers, and the denominators have at least one zero.

\begin{theorem}\label{thm:FerrersRatio}
	Let $\mu$, $\nu \in \mathbb{R}$. If $N_1 \geq 1$ and $N_2 \geq 2$ (or $N_2 \geq 1$ and $N_1 \geq 2$), then the function $x \mapsto y_1(x)/y_2(x)$ (or $x \mapsto  y_2(x)/y_1(x)$) is strictly decreasing (or increasing) between two consecutive zeros of $y_2(x)$ (or $y_1(x)$).
\end{theorem}

Note that the volume constraint in Condition 5.5 always implies $r_1 + r_2 < \pi/\sqrt{k}$. In \cite{SKAK24}, the authors analytically verified Condition 5.5 for the case $n = \beta/k$ by showing that the $l$ constraints imply $r_1 + r_2 = \pi/\sqrt{k}$, contradicting the inequality $r_1 + r_2 < \pi/\sqrt{k}$. In the case when $n > \beta/k$, their numerical evaluations (\cite[p. 3317]{SKAK24}) suggest that the $l$ constraints force $r_1 + r_2 > \pi/\sqrt{k}$, which yields an indication of contradiction to the volume constraints. To handle this case with our technique, we need the ratio $y_1/y_2$ to be monotonic on $(-1, x_{2, N_2})$, with a specific restriction on $\mu$, which imposes further restrictions on $n$, $\beta$ and $k$. However, we can show that $n>\beta/k$ is impossible in the setting of Condition 5.5, by using the technique of Lemma \ref{lemma:fDashNozero}. Same numerical tests in case of $n < \beta/k$ suggest that the $l$ constraints also implies  $r_1 + r_2 < \pi/\sqrt{k}$. In this paper, we verify Condition 5.5 whenever $n < \beta/k$, the case where the numerical tests do not indicate a contradiction.

Notice that, if $n < \beta/k$ and  $\nu=\gamma_{\beta,k,n},$ then $0<N_1<N_2$. Now, we observe the distribution of zeros directly from Lemma \ref{lemma:interlace} as follows.

\begin{remark}\label{remark1}
	Suppose that $y_1$ and $y_2$ have unequal strictly positive numbers of zeros $N_1$ and $N_2$; then the arrangement of zeros is as follows.
	\begin{enumerate}
		\item[\bf a.] If $N_1 < N_2$, then $-1 < x_{2,1} < x_{1,1} < x_{2,2}<\ldots <x_{1, N_1} < x_{2, N_2}<1.$
		\item[\bf b.] If $N_2 < N_1$, then $-1 < x_{1,1} < x_{2,1} < x_{1,2}<\ldots <x_{2, N_2} < x_{1, N_1}<1.$
	\end{enumerate}
\end{remark}

From the above remark it is clear that if $\nu= \gamma_{\beta,k,n},$ then $0< \tilde{R}=r_{2, N_2}< r_{1, N_1}$. Now, from Lemma \ref{lemma:fDashNozero} we have  that $\tilde{R}_2< r_{1, N_1}$. Therefore we have $(0, r_{2, N_2})= (0, \tilde{R})$ and $(\tilde{R}, \tilde{R}_2) \subset (r_{2, N_2}, r_{1,N_1})$. Moreover, from the next theorem we obtain the monotonicity and the sign of $l$ in these intervals.

\begin{theorem}\label{thm:l-monotonicityk>0}
	If $n < \beta/k$, that is $0<N_1<N_2$, then the function $l$ is positive and strictly increasing on $(0, r_{2, N_2})$ and strictly decreasing on $(r_{2, N_2}, r_{1, N_1})$.
\end{theorem}
	
\subsection{Verification of Condition 5.5 for $k>0$}

From Theorem \ref{thm:l-monotonicityk>0} we claim that, $l$ is positive and increasing on $(0, \tilde{R})$ and positive and decreasing on $(\tilde{R}, \tilde{R}_2)$. Notice that $l(0)=0=l(r_{1,N_1})$ and  $\tilde{R}_2 < r_{1,N_1}$. Hence, there exists a point $\tilde{R}_2^* \in (0, \tilde{R})$ such that $l(\tilde{R}_2^{*})=l(\tilde{R}_2)$. Therefore for each $r_1 \in (\tilde{R}_2^*, \tilde{R})$ there exists a unique $r_2 \in (\tilde{R}, \tilde{R}_2)$ such that $l(r_1)=l(r_2)$. We do not work with the interval $(0, \tilde{R}_2^*)$ as for every point $r_1$ in this interval there are no points $r_2 \in (\tilde{R}, \tilde{R}_2)$ such that $l$ constraints holds. Hence, in this interval, Condition 5.5 is trivially true. Also notice that the rate of change of $r_2$ with $r_1$ is negative, that is ${\rm d}r_2/{\rm d}r_1= {l'(r_1)}/{l'(r_2)}<0$. Therefore as $r_1$ moves from $\tilde{R}_2^*$ to $\tilde{R}$, its corresponding point $r_2$ moves from $\tilde{R}_2$ to $\tilde{R}$. There is a special radius $R_{1/2}$ (see \cite[p. 3313]{SKAK24}) for each $R \in (\tilde{R}, \tilde{R}_2)$ such that $V(R_{1/2})= V(R)/2$. Notice that $V$ is an increasing function. Hence, we have $V(R)<V(\tilde{R}_2)=2V(\tilde{R})$. This gives $V(R_{1/2})= V(R)/2<V(\tilde{R})$. Therefore we obtain  $R_{1/2}<\tilde{R}$. We now would like to point out that $\tilde{R}_2^*<R_{1/2}$. Indeed, if for the sake of contradiction we would assume that $R_{1/2}\leq\tilde{R}_2^*$, then we would obtain that this cannot happen since in this case the $l$ constrains fails as we observed above. Now, we define two new points $R_{1/2}^{*}$ in $(\tilde{R}, \tilde{R}_2)$ and $R^{*}$ in $(0, \tilde{R})$ such that $l(R_{1/2}^{*})=l(R_{1/2})$ and $l(R^{*})=l(R)$. We are left with two possible arrangements of these numbers:
\begin{equation}\label{case1}
	\tilde{R}_2^{*}<R_{1/2}<R^{*}<\tilde{R}< R< R^{*}_{1/2}< \tilde{R}_2
\end{equation}
or
\begin{equation}\label{case2}
	\tilde{R}_2^{*}<R^{*}<R_{1/2}<\tilde{R}< R^{*}_{1/2}< R< \tilde{R}_2.
\end{equation}
Let us deal with this case by case. Let $r_1<\tilde{R}< r_2$ such that $l(r_1)=l(r_2)$.

{\bf Case 1}. We first consider the case when \eqref{case1} holds. If $r_1< R^{*},$ then clearly $l(r_2)=l(r_1)< l(R^{*})= l(R)$ as $l$ is an increasing function on $(0, \tilde{R})$. Since $l$ is a decreasing function on $(\tilde{R}, \tilde{R}_2)$ we have that $R<r_2$. Now, since $V$ is an increasing function, we obtain that $V(R)< V(r_2)$ and hence $V(r_1)+V(r_2)> V(R)$. If $r_1> R^{*},$ then evidently $R_{1/2}< r_1<r_2$. Hence $V(r_1)+V(r_2)> 2 V(R_{1/2})=V(R)$. In this arrangement, we proved that Condition 5.5 is valid and hence the statement in Conjecture 5.2 of \cite{SKAK24} is true.
\begin{figure}[ht]
	\centering
	\begin{tikzpicture}[xscale=0.85, yscale=0.55, every node/.style={font=\scriptsize}]
		\draw[-latex, thick] (0,0) -- (10.5,0) node[right] {$r$};
		\draw[-latex, thick] (0,0) -- (0,5.2) node[above] {$l(r)$};
		\node[below left] at (0,0) {$0$};
		\draw[ultra thick, black, name path=curveL] (0,0) to[out=15, in=260] (3.35, 5.0);
		\draw[ultra thick, black, name path=curveR] (3.65, 5.0) to[out=280, in=165] (9.5, 0);
		\draw[dashed, black, ultra thick] (3.5,0) -- (3.5,5.2);
		\node[above] at (3.5, 5.2) {Pole $\tilde{R}$};
		\node[below] at (3.5,-0.1) {$\tilde{R}$};
		\path[name path=line1] (0,1.5) -- (9.5,1.5);
		\path[name path=line2] (0,2.3) -- (9.5,2.3);
		\path[name path=line3] (0,3.8) -- (9.5,3.8);
		\path[name intersections={of=curveL and line1, by=L1}];
		\path[name intersections={of=curveR and line1, by=R1}];
		\draw[dashed, gray, thick] (L1) -- (R1);
		\draw[dashed, gray] (L1) -- (L1 |- 0,0);
		\draw[dashed, gray] (R1) -- (R1 |- 0,0);
		\filldraw[gray] (L1) circle (1.5pt);
		\filldraw[gray] (R1) circle (1.5pt);
		\filldraw[gray] (L1 |- 0,0) circle (1.5pt);
		\filldraw[gray] (R1 |- 0,0) circle (1.5pt);
		\node[below] at (L1 |- 0,0) {$\tilde{R}_2^*$};
		\node[below] at (R1 |- 0,0) {$\tilde{R}_2$};
		\path[name intersections={of=curveL and line2, by=L2}];
		\path[name intersections={of=curveR and line2, by=R2}];
		\draw[dashed, black!60!white, thick] (L2) -- (R2);
		\draw[dashed, black!60!white] (L2) -- (L2 |- 0,0);
		\draw[dashed, black!60!white] (R2) -- (R2 |- 0,0);
		\filldraw[black!60!white] (L2) circle (1.5pt);
		\filldraw[black!60!white] (R2) circle (1.5pt);
		\filldraw[black!60!white] (L2 |- 0,0) circle (1.5pt);
		\filldraw[black!60!white] (R2 |- 0,0) circle (1.5pt);
		\node[above=2pt] at (L2 |- 0,0) {$R_{1/2}$};
		\node[above=2pt] at (R2 |- 0,0) {$R^*_{1/2}$};
		\path[name intersections={of=curveL and line3, by=L3}];
		\path[name intersections={of=curveR and line3, by=R3}];
		\draw[dash dot, black, thick] (L3) -- (R3);
		\draw[dashed, black] (L3) -- (L3 |- 0,0);
		\draw[dashed, black] (R3) -- (R3 |- 0,0);
		\filldraw[black] (L3) circle (1.5pt);
		\filldraw[black] (R3) circle (1.5pt);
		\filldraw[black] (L3 |- 0,0) circle (1.5pt);
		\filldraw[black] (R3 |- 0,0) circle (1.5pt);
		\node[below left, inner sep=1pt] at (L3 |- 0,0) {$R^*$};
		\node[below right, inner sep=1pt] at (R3 |- 0,0) {$R$};
		\node[below] at (9.5,-0.1) {$r_{1,N_1}$};
		\end{tikzpicture}
		\caption{Illustration of the arrangement in Case 1.}
		\label{fig:radii}
	\end{figure}
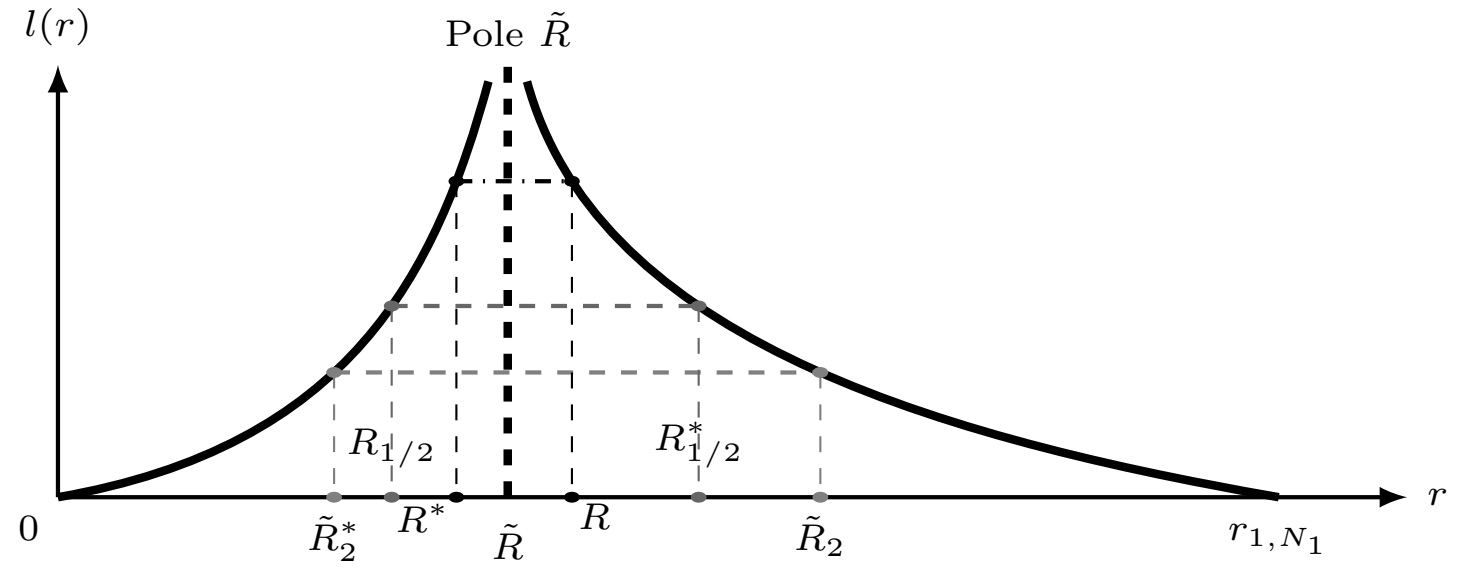
	
	{\bf Case 2}. Now, we consider the case when \eqref{case2} holds. In this case for every $r_1<R^{*}$ we have $r_2> R$ and hence $V(r_1)+V(r_2)>V(R)$. If we suppose that $R_{1/2}<r_1,$ then $V(r_1)+V(r_2)> 2 V(R_{1/2})=V(R)$. In case of $r_1 \in (R^{*}, R_{1/2})$ which strictly maps to $r_2 \in (R_{1/2}^{*}, R)$, we could not obtain similar contradiction involving volume constraints. Thus, we have rigorously verified Condition 5.5 for the disjoint sub-domains $r_1 \in (0, R^{*}] \cup [R_{1/2}, \tilde{R})$, establishing a partial verification of the conjecture.

\begin{remark}\label{remark:case2Improvement}
Now, we are going to construct a strictly decreasing sequence $\{c_k\}$ such that $V(r_1)+V(r_1^*)>V(R)$ for every $r_1 \in (c_k, R_{1/2})$ for every $k$. Let $r \in (0, \tilde{R})$ and $r^* \in (\tilde{R}, \tilde{R}_2)$ such that $ l(r) = l(r^*)$. We define $G: (0, \tilde{R}) \to \mathbb{R}$ by
$$G(r) = V(r) + V(r^*) - V(R).$$
Since both functions $V$ and $r \mapsto r^*$ are continuous, $G$ is also continuous. Note also that $G(R_{1/2})>0$. Therefore there exists some $\delta > 0$ such that we have $ G(r) > 0$ for every  $r \in [R_{1/2} - \delta, R_{1/2}]$. We choose our initial starting point $c_0$ to lie within this interval, that is $ c_0 \in [R_{1/2} - \delta, R_{1/2}]$. This guarantees that $G(c_0) > 0$. We define the sequence $\{c_k\}_{k \ge 0}$ starting at $c_0$ by $ V(c_{k+1}) = V(R) - V(c_k^*)$. We prove that the sequence $\{c_k\}$ is strictly monotonically decreasing, by using mathematical induction. First, we show that $c_1 < c_0$. By definition, we have
$$ V(c_1) = V(R) - V(c_0^*). $$
Since $G(c_0) > 0$, we have $V(c_0^*) > V(R) - V(c_0)$. Therefore we obtain
$$ V(c_1) = V(R) - V(c_0^*) < V(R) - \left(V(R) - V(c_0)\right) = V(c_0). $$
Since $V$ is strictly increasing, $V(c_1) < V(c_0)$ implies $ c_1 < c_0 $. Now assume that $c_k < c_{k-1}$ for some $k \ge 1$. Since $l$ is strictly increasing on $(0, \tilde{R})$, we have $l(c_k) < l(c_{k-1}).$ This implies $ l(c_k^*) = l(c_k) < l(c_{k-1}) = l(c_{k-1}^*)$. Since $l(r)$ is strictly decreasing on $(\tilde{R}, \tilde{R}_2)$,  we have $ c_k^* > c_{k-1}^*  $. Therefore, $ V(c_k^*) > V(c_{k-1}^*)$. This yields
$$ V(R) - V(c_k^*) < V(R) - V(c_{k-1}^*).$$
By definition, the above inequality is exactly $ V(c_{k+1}) < V(c_k)$. Since $V(r)$ is strictly increasing, this proves $ c_{k+1} < c_k$. By mathematical induction, the sequence $\{c_k\}_{k \ge 0}$ is strictly monotonically decreasing
$$c_0 > c_1 > c_2 > \ldots > c_k > c_{k+1} > {\ldots}. $$
		
Again, by induction, we show $G(r) > 0$ for all $r > c_k$. For any $r \ge R_{1/2}$, we already know that $G(r) > 0$. Moreover, the same is true for any $r \in (c_0, R_{1/2})\subset [R_{1/2} - \delta, R_{1/2}],$ since we have guaranteed positivity by our choice of $\delta$. Thus, we have that $G(r)>0$ for all $r>c_0.$ Now, assume that $G(r) > 0$ for all $r > c_k$ for some $k$. Since the induction hypothesis already guarantees that $G(r) > 0$ for all $r > c_k$, we only need to prove positivity on the remaining gap $r \in (c_{k+1}, c_k]$. Since $r > c_{k+1}$, we have
$$ V(r) > V(c_{k+1}) = V(R) - V(c_k^*). $$
On the other hand, the condition $r \le c_k$ implies that $l(r) \le l(c_k) = l(c_k^*).$ Since $l$ is strictly decreasing on the right branch, we have $ r^* \ge c_k^*,$ which implies $V(r^*) \ge V(c_k^*)$. Adding these two volume inequalities together, we obtain
$$ V(r) + V(r^*) > \left(V(R) - V(c_k^*)\right) + V(c_k^*) = V(R) $$
which simplifies directly to
$$ G(r) = V(r) + V(r^*) - V(R) > 0. $$
By mathematical induction, $G(r) > 0$  for all $r > c_k$ for every $k \ge 0$.
Since $c_k$ decreases strictly, it can cross the threshold $R^*$ in a finite number of steps (i.e., there exists some finite index $m \ge 1$ such that $c_m \le R^*$). If this occurs, then we are done. However, if the sequence ${c_k}$  never crosses, then it is bounded below by $R^*,$ and hence it is convergent. If the sequence ${c_k}$ were to converge to a limit $L$ strictly greater than $R^*$, the sequence would indeed get trapped at $L.$
\end{remark}
	
\subsection{Monotocity of the ratio $y_3/y_4$ involving conical functions}

Now we obtain the monotonicity below through a technique identical to the one used in Theorem \ref{thm:FerrersRatio}. For this reason we omit the proof of Theorem \ref{thm:conicalRatio}.
\begin{theorem}\label{thm:conicalRatio}
	Let $\tau \neq 0$ and $\mu \in \mathbb{R}$.
	The function $ x \mapsto y_3(x)/y_4(x)$ (or $ x \mapsto  y_4(x)/y_3(x)$) is strictly increasing (or decreasing) between two consecutive zeros of $y_4(x)$ (or $y_3(x)$).
\end{theorem}
	
Before proceeding to the verification of Condition 5.5, we obtain the distribution of zeros of conical functions. This distribution gives the ordering of the singularities of the $l$ function in the $r$ domain.
\begin{theorem} \label{thm:conicalZeros}The distributions of zeros is as follows:
	\begin{enumerate}
		\item[\bf a.] If $\mu \geq 0,$ then $x_{4,1}<x_{3,1}.$
		\item[\bf b.] If $\mu < 0,$ then $x_{3,1}<x_{4,1}$.
	\end{enumerate}
\end{theorem}

\subsection{Verification of Condition 5.5 for $k<0$}

Let $r_{3,j}$ be the corresponding zero in $r$ domain for the zero $x_{3,j}$ in $x$ domain for $j\geq1.$ Similarly, let $r_{4,j}$ be the corresponding zero to $x_{4,j}$ for $j\geq1.$ Since $\sinh(\sqrt{-k}r)>0$ for all $ r \in (0, \infty)$, from the definitions of $f$ and $l$ it is clear that the zeros $r_{4,j}$ of $f$ are the singularities of $l$ for $j\geq1.$ Notice that if $x \in \left(1,\infty \right)$ then $ r \in \left(0, \infty \right)$ whenever $x=\cosh(\sqrt{-k}r)$, and the hyperbolic cosine function is increasing on $\left[0, \infty\right)$. Hence for every $j$ we have $r_{3,j}< r_{3,j+1}$ and $r_{4,j}< r_{4,j+1}$.
Therefore in case of the profile function $f$, the first zero $\tilde{R}$ is $r_{4,1}$.
	
Now, from part {\bf a} of Theorem \ref{thm:conicalZeros} we clearly have that $r_{4,1}<r_{3,1}<r_{4,2},$ and from Lemma \ref{lemma:fDashNozero} we obtain that $\tilde{R}_2< r_{3,1}$. Thus, we arrive to the following Theorem \ref{thm:l-monotonicityk<0}, in whose proof we can use a very similar technique as we used in the proof of Theorem \ref{thm:l-monotonicityk>0}. Because of this, we will omit the proof of Theorem \ref{thm:l-monotonicityk<0}.

\begin{theorem}\label{thm:l-monotonicityk<0}
The function $l$ is positive and strictly increasing on $(0, r_{4, 1})$ and strictly decreasing on $(r_{4, 1}, r_{3, 1})$.
\end{theorem}
Now the rest of the verification is similar to the case of $k>0$.
	
\begin{remark}
During the process of finalising this paper, we found that \cite[eq. (25)]{SKAK24} is not correct, and the correct version is as follows:
\begin{equation}
	l(r) = \begin{cases}
	\left( \frac{y_1(x)}{y_2(x)} \right)^2, & \text{for } x = \cos(\sqrt{k}r) \in (-1,1),\ k > 0, \\[1.5ex]
	\left( \frac{y_3(x)}{y_4(x)} \right)^2, & \text{for } x = \cosh(\sqrt{-k}r) \in (1, \infty),\ k < 0.
	\end{cases}
\end{equation}
This can be verified, directly from the differential equation \cite[eq. 22]{SKAK24} of the profile function $f,$ by showing that in both cases of $k\lessgtr0$ we have $j'(r)/h(r) = -\frac{1}{\beta} (f'(r)/f(r))^2.$  It is worth mentioning that, unfortunately, this correction does not help us to fill the gap in verifying the corrected Condition 5.5, which would settle \cite[Conjecture 5.2]{SKAK24}. Observe that, if $k > 0,$ then
$$l'(r)= -2\sqrt{k}\sqrt{1-x^2} \left(\frac{y_1(x)}{y_2(x)}\right) \frac{{\rm d}}{{\rm d}x} \left(\frac{y_1(x)}{y_2(x)}\right)$$
for all $x = \cos(\sqrt{k}r) \in (-1,1)$ and if $k<0,$ then
$$ l'(r)= 2\sqrt{-k}\sqrt{x^2-1} \left(\frac{y_3(x)}{y_4(x)}\right) \frac{{\rm d}}{{\rm d}x} \left(\frac{y_3(x)}{y_4(x)}\right)$$
for all $x = \cosh(\sqrt{-k}r) \in (1, \infty)$. Hence, with the above corrections the statements of Theorems \ref{thm:l-monotonicityk>0} and \ref{thm:l-monotonicityk<0} are still valid and the modified definition of $l(r)$ yields the same conclusions as those presented in our paper.
\end{remark}
	
\section{Proofs of the Main Results}
	
We use the following lemma to observe that $y_3$ and $y_4$ have infinitely many zeros.

\begin{lemma}\label{lemma:infintezeros}
If $\tau \neq 0$ then $y_3$ and $y_4$ have infinitely many zeros on $(1, \infty)$.
\end{lemma}

\begin{proof}[\bf Proof]
By setting $x= \cosh(r)$ and $\nu=\degree$ we transform the equation \eqref{eq:legendreequation} into the differential equation
$$\omega''(r) + (\coth r) \omega'(r) + \left(\tau^2+\frac14 - \frac{\mu^2}{\sinh^2 r}\right)\omega(r) = 0$$
for the function $r \mapsto \omega(r)=y(\cosh(r))$ on $(0, \infty)$. Moreover, observe that with
$$\omega(r)= u(r)\exp\left(-\frac{1}{2}\int\coth(r) {\rm d}r\right)= \frac{u(r)}{\sqrt{\sinh (r)}}$$
we transform the above equation into $$u''(r) + \left( \tau^2 + \frac{\frac14 - \mu^2}{\sinh^2 r} \right) u(r) = 0,$$ where $r \in (0, \infty)$. Thus, if we apply Leighton's Lemma \cite[Theorem 31.7]{AgarwalORegan} we conclude that every solution of the above differential equation has infinitely many zeros on $(0,\infty),$ which in turn implies that indeed the functions $y_3$ and $y_4$ have infinitely many zeros on $(1, \infty).$
\end{proof}
	
The next lemma verifies the validity of \eqref{eqtaudef}.

\begin{lemma}\label{lemma:Mckean}
If $k<0,$ then $\gamma_{\beta,k,n}$ is of the form $\degree$ with $\tau=\tau_{\beta,k,n}$ defined in \eqref{eqtaudef}.
\end{lemma}

\begin{proof}[\bf Proof]
Let $B(r)$ denote a geodesic ball of radius $r$, and let $\lambda_1(B(r))$ denote the first eigenvalue of \cite[eq. 12]{SKAK24} on $B(r)$. Since, for every $\epsilon>0$ we have that $B(\tilde{R}+\epsilon)\setminus\overline{B(\tilde{R})}$ is non empty and open, it follows that $\beta=\lambda_1(B(\tilde{R}))>\lambda_1(B(\tilde{R}+\epsilon))$ (see \cite[p. 18]{Chavel}). On the other hand, we know that $B(\tilde{R}+\epsilon)$ is a regular domain, hence a normal domain. Thus, we have that $M^{n}_{k}$ is simply connected. Therefore as a conseuence of McKean's theorem \cite[p. 70]{Chavel} we obtain that $-\beta=-\lambda_1(B(\tilde{R}))<-\lambda_1(B(\tilde{R}+\epsilon))\leq\frac{(n-1)^{2}k}{4}<0,$ by following the definitions of regular and normal domains from \cite{Chavel}.
\end{proof}
	
Now, we are going to present Lemmas \ref{lemma:simplezeros} - \ref{lemma:interlace}, together with proofs for completeness. We use these lemmas in the proofs of Theorem \ref{thm:FerrersRatio} and Theorem \ref{thm:conicalRatio}.

\begin{lemma}\label{lemma:simplezeros}
Every zero of $y_1$ and $y_2$ in the open interval $(-1, 1)$, and of $y_3$ and $y_4$ in the open interval $(1, \infty)$, is simple.
\end{lemma}

\begin{proof}[\bf Proof]
Since $y_1, y_2, y_3,$ and $y_4$ are the non-trivial solutions of \eqref{eq:legendreequation}, we conclude the proof by using the classical contradiction method used for the solutions of second order homogeneous differential equations (see for example \cite[p. 139]{OlvASF97}).
\end{proof}

\begin{lemma}\label{lemma:NoCommonZeros}
The functions $y_1$ and $y_2$ cannot have a common zero in $(-1,1)$, and the functions $y_3$ and $y_4$ cannot have a common zero in $(1,\infty)$.
\end{lemma}

\begin{proof}[\bf Proof]
Suppose, for contradiction, that there exists $c \in (-1, 1)$ such that $y_1(c) = y_2(c) = 0$. In view of equations \eqref{eq:y1coupledpart} and \eqref{eq:y2coupledpart}, we arrive at $(1-c^2) y_1'(c)=0$ and $(1-c^2) y_2'(c)=0$. These in turn imply that $y_1'(c) = y_2'(c) = 0,$ which contradicts the fact that $y_1$ and $y_2$ has simple zeros (see Lemma \ref{lemma:simplezeros}). Hence, $y_1$ and $y_2$ cannot share a common zero in $(-1,1)$. By similar arguments using the equation \eqref{eq:y3coupledpart} and \eqref{eq:y4coupledpart}, we conclude that  $y_3$ and $y_4$ cannot share a common zero in $(1,\infty)$.
\end{proof}

\begin{lemma}\label{lemma:interlace}
The zeros of $y_1$ and $y_2$ are strictly interlaced in the interval $(-1, 1)$, and similarly, the zeros of $y_3$ and $y_4$ are strictly interlaced in the interval $(1, \infty)$.
\end{lemma}

\begin{proof}[\bf Proof]
Let $\mu$, $\nu$, $\tau$ be any real numbers. Let $a<b$ be two consecutive zeros of $y_{2}$. Then clearly $y_2'(a)y_2'(b)<0$. Equation \eqref{eq:y2coupledpart} gives $$y_2'(a)y_2'(b)=\frac{1}{\sqrt{1-a^2}\sqrt{1-b^2}}(\mu+\nu)^{2}(\nu+1-\mu)^{2}y_1(a)y_1(b).$$ Now, we have $y_1(a)y_1(b)<0$. Continuity of $y_1$ ensures that there is at least one zero of $y_1$ in $(a,b)$. By similar arguments, in view of the equation \eqref{eq:y1coupledpart}, we obtain that at least one zero of $y_2$ is between two consecutive zeros of $y_1$. Again, with similar arguments for $y_1$ and $y_2$, and using \eqref{eq:y3coupledpart} and \eqref{eq:y4coupledpart}, we conclude the remaining part of the lemma.
\end{proof}
	
We shall use the next lemma in the proof of Lemma \ref{lemma:keyInequality}, which we require in the proof of Theorem \ref{thm:FerrersRatio}.

\begin{lemma}\label{lemma:NistMissing}
Let $\alpha >0$ and $\nu \in \mathbb{R}$ be order and degree of the Ferrers function $x \mapsto \mathsf{P}^{\alpha}_{\nu}(x)$ on $(-1, 1)$. If the integer part of the order $ \lfloor\alpha\rfloor$, the integer part of the degree $\lfloor \nu \rfloor$, the fractional parts $\delta_{\nu}=\nu -\lfloor \nu \rfloor$ and $\delta_{\alpha}= \alpha-\lfloor\alpha\rfloor$ satisfy $\lfloor \alpha \rfloor = \lfloor \nu \rfloor$ and $0<\delta_{\nu} \leq \delta_{\alpha}$, then number of zeros of $x \mapsto \mathsf{P}^{\alpha}_{\nu}(x)$ is $1$.
\end{lemma}

\begin{proof}[\bf Proof]
From the series representation \cite[eq. 14.3.1]{Nist} of $\mathsf{P}_\nu^\alpha(x)$ it is clear that every zero of it on $(-1, 1)$ is in one-to-one correspondence with the zeros of
${}_2F_1\left(\nu+1, -\nu, 1-\alpha; \frac{1-x}{2}\right)={}_2F_1(a, b, c; t)$, where $a=\nu+1$, $b=-\nu$ and $c=1-\alpha$.
Throughout the proof we assume $x \in (-1,1)$ and $t=\frac{1-x}{2}$. If $N$ is the common integer part $N=\lfloor\alpha\rfloor=\lfloor \nu \rfloor,$ then $$f^{(N)}(t) = \frac{(a)_N (b)_N}{(c)_N} {}_2F_1(a+N, b+N; c+N; t)= \frac{(a)_N (b)_N}{(c)_N}\left(1 + \sum_{k\geq1}c_j t^k\right),$$ where $c_k = \frac{(a+N)_k (b+N)_k}{(c+N)_k k!}<0$. Notice that, $f^{(N)}(0)=\frac{(a)_N (b)_N}{(c)_N}>0 $. Therefore $f^{(N)}$  is strictly decreasing on $(0, 1)$. Since $f^{(N)}(0) > 0$ and  $c+N - (a+N) - (b+N)< 0$ the series $1 - \sum_{k\geq1}|c_k| t^k$ does not converge at $t = 1$. Therefore $\lim_{t \nearrow 1}f^{(N)}(t)=-\infty$. Hence  $f^{(N)}(t)$ has exactly 1 zero in the interval $(0, 1)$.

On the other hand, observe that for $0\leq j \leq N$ we have $f^{(j)}(0) >0 $. Now, we prove that if for some $j\in\{1,\ldots,N\}$ we have that $f^{(j)}$ has only one zero in $(0, 1)$ then $f^{(j-1)}$ has only one zero in $(0, 1)$. Let $t_0$ be one and only zero of $f^{(j)}$ in $(0, 1)$, for some $j\in\{1,\ldots,N\}$. Since, $f^{(j)}(0)>0$, we have $f^{(j)}(t)>0$ on $[0, t_0)$ and $f^{(j)}(t)<0$ on $( t_0, 1)$. Therefore $f^{(j-1)}$ is strictly increasing on $[0, t_0)$ and decreasing on $(t_0, 1)$. Now, we have $0<f^{(j-1)}(0)<f^{(j-1)}(t_0)$ and due to the divergence of the hypergeometric part, we observe that $\lim_{t\nearrow 1}f^{(j-1)}(t)=-\infty$. This implies $f^{(j-1)}$ has exactly one zero in $(t_0, 1)$ and hence in $(0,1)$. Consequently $f^{(0)}=f$ has only one zero in $(0, 1)$.
\end{proof}
	
\begin{lemma}\label{lemma:keyInequality}
If $N_2 \geq 2 $, then $(\mu+\nu)(\nu-\mu+1) > 0$.
\end{lemma}

\begin{proof}[\bf Proof]
Let $\lceil x \rceil$ be the least integer greater than or equal to $x \in  \mathbb{R}$. From the number of zeros formula in \cite[p. 368]{Nist} and Lemma \ref{lemma:NistMissing} we observe that the number $N_2$ is given by $\max(\lceil \nu - |1- \mu| \rceil, 0)+1$ or $\max(\lceil \nu - |1- \mu| \rceil, 0)$. Now $N_2 \geq 2$ implies $ \max(\lceil \nu - |1- \mu| \rceil, 0) \geq 1$. Therefore $\lceil \nu - |1- \mu| \rceil \geq 1$ gives $\nu - |1- \mu|>0$. Hence $\mu-\nu-1< 0 < \mu+\nu$.
\end{proof}

\begin{proof}[\bf Proof of Theorem \ref{thm:FerrersRatio}]
We first define a mapping $R: (a,b) \to \mathbb{R}$ by
$$R(x) = \frac{y_{1}(x)}{y_{2}(x)} \quad \text{for every } x \in (a,b),$$
where $a$ and $b$ are two consecutive zeros of $y_{2}(x)$ on $(-1, 1)$. Since $y_1$ and $y_2$ are differentiable on $(a,b)$, $R$ is also differentiable on $(a,b)$ and the derivative of $R$ at any $x$ in $(a,b)$ is clearly $(y'_{1}(x)y_{2}(x)-y'_{2}(x)y_{1}(x))/(y_{2}(x))^{2}$. Now, define a continuously  differentiable Wronskian type function  $W:[a,b] \to \mathbb{R}$ by
$$W(x)=y'_{1}(x)y_{2}(x)-y'_{2}(x)y_{1}(x) \quad \text{for every } x \in [a,b].$$
Observe that, to show that $R$ is decreasing in $(a,b)$, it is sufficient to show that $W(x)<0$ for every $x \in (a,b)$. From definition of $W$, Lemma \ref{lemma:keyInequality}, and equations \eqref{eq:y1coupledpart} and \eqref{eq:y2coupledpart} we obtain that
$$W(a) = -y_{2}'(a)y_{1}(a)=\frac{-(\mu+\nu)(\nu-\mu+1)}{\sqrt{1-a^2}}(y_1(a))^2<0$$
and
$$W(b)=-y_{2}'(b)y_{1}(b)=\frac{-(\mu+\nu)(\nu-\mu+1)}{\sqrt{1-b^2}}(y_1(b))^2 <0.$$
Therefore if $W$ has no zeros in $(a,b)$ then $W(x)<0$ for every $x \in (a,b)$. Next, we show that $W$ does not vanish in $(a,b)$, by considering two cases $\mu =1/2$ and $\mu \neq 1/2$.
	
First, if $\mu =1/2$, then both $y_1$ and $y_2$ are solution of the same differential equation
$$(1-x^2)y''(x)-2xy'(x)+ \left[\nu(\nu+1)-\frac{1}{4 (1-x^2)}\right]y(x)=0.$$
Hence $W(x)$ is $-1$ times of the Wronskian of $y_1$ and $y_2$, where in this case $y_{1}(x)=\mathsf{P}^{-1/2}_{\nu}(x)$ and $y_2(x)= \mathsf{P}^{1/2}_{\nu}(x).$ From \cite[(14.2.7)]{Nist}, we have that the Wronskian of $\mathsf{P}^{-1/2}_{\nu}(x)$ and $\mathsf{P}^{1/2}_{\nu}(x)$ becomes $\frac{2}{\pi (1-x^2)}>0$ for every $x \in (a,b) \subset (-1,1)$.
Therefore if $\mu=1/2$ then $W(x)<0$ for every $x \in (a,b)$.
	
Next, for the case of $\mu \neq 1/2$, we prove that there are no zeros of $W$ in $(a,b)$ through the contradiction method. Assume for contradiction that there is a point $d$ in $(a,b)$ such that $W(d)=0$. We claim that if such $d$ exists, then $W'(d) \neq 0$. To prove the claim, we first construct a differential equation satisfied by $W$. Since $y_1$ and $y_2$ are solutions of \eqref{eq:legendreequation} on $(-1,1)$, we obtain that
$$(1 - x^2)y_{1}''(x)y_{2}(x) - 2xy_{1}'(x) y_{2}(x) + \left[ \nu (\nu + 1) - \frac{\mu^2}{1 - x^2} \right]y_{1}(x) y_{2}(x) = 0,$$
$$(1 - x^2)y_{2}''(x) y_{1}(x) - 2xy_{2}'(x) y_{1}(x) + \left[ \nu (\nu + 1) - \frac{(1-\mu)^2}{1 - x^2} \right]y_{2}(x) y_{1}(x) = 0,$$
where $ x \in [a,b]$, $\mu \neq 1/2$ and $\nu \in \mathbb{R}$. Combining the above equations, we obtain the differential equation
\begin{equation}\label{W-DE}
(1-x^2)\zeta'(x)-2x \zeta(x)+\left(\frac{1-2\mu}{1-x^2}\right)y_{1}(x)y_{2}(x)=0
\end{equation}
satisfied by $W$ on $[a,b]$, where $\mu \neq 1/2$ and $\nu \in \mathbb{R}$. Therefore, from \eqref{W-DE} we have
$$(1-d^2)W'(d)=-\left(\frac{1-2\mu}{1-d^2}\right)y_{1}(d)y_{2}(d).$$ Since $\mu \neq 1/2$ and $d \in (a,b)$, we have $(1-2\mu) \neq 0$ and $y_2(d) \neq 0$. We claim that $y_1(d) \neq 0$. Suppose that $y_1(d)=0$, then by definition of $W$ and \eqref{eq:y1coupledpart} we have
$$W(d)= y_1'(d)y_2(d)= \frac{-(y_2(d))^{2}}{\sqrt{1-d^2}}<0.$$
However $W(d) = 0$. Therefore, we attain $0<0$, which is absurd. Hence $y_1(d) \neq 0$. This implies $W'(d) \neq 0$.
		
We also claim that $d \neq 0$. For the sake of contradiction we suppose that $d=0,$ that is $W(0)=0,$ which clearly implies that $y_1(d)\neq0$ and $y_2(d)\neq 0$. However, from \eqref{eq:y1coupledpart} and \eqref{eq:y2coupledpart} we have $W(0)= -((y_2(0))^2 + (\mu+\nu)(\nu-\mu+1)(y_1(0))^2) \neq 0$. This is again a contradiction. Therefore $d \neq 0$.
	
Since $W$ is a differentiable function on $(a,b)$ with simple zeros, every zero is isolated. Then $W$ must have only finitely many zeros in $(a,b)$. Otherwise there exist a convergent sequence $\{x_n\}_{n \in \mathbb{N}}$ of zeros of $W$ in $[a,b]$. However, continuity of $W$ forces the limit point $x_0$ of such a sequence to be a non-isolated zero of $W$. Therefore, no interior points of $[a,b]$ are limit points. This forces either $a$ or $b$ to be the limit point $x_0$. Hence $W(a)=0$ or $W(b)=0$. This is in contradiction with $W(a)<0$ and $W(b)<0$. Let $d_1$ be the first zero of  $W$ in $(a,b)$. Since $W(a)<0$,  we have $W'(d_1)>0$. Notice that $W(b)<0$. This forces the existence of consecutive zero $d_2$ to $d_1$ such that $d_1 <d_2$ and  $W'(d_2)<0$. Now we show that, if $d$ is a zero of $W$, then $d$ and $W'(d)$ have the same sign. For this, by using \eqref{eq:y1coupledpart} and \eqref{eq:y2coupledpart} we obtain
\begin{align*}
(1-d^2)W'(d)&=-\left(\frac{1-2\mu}{1-d^2}\right)y_{1}(d)y_{2}(d) \\
&=\frac{(\mu+\nu)(\nu-\mu+1)\sqrt{1-d^2}(y_1(d))^{2}+\sqrt{1-d^2}(y_2(d))^{2}}{d(1-d^2) }\neq 0.
\end{align*}
From the above equation, Lemma \ref{lemma:keyInequality} and the fact that $d \in (a,b) \subset (-1,1)$, we conclude that both $d$ and $W'(d)$ have the same sign. Since $W'(d_1)>0$ and $W'(d_2)<0$ we have $d_2<0<d_1$. This is a contradiction to our construction of $d_1 <d_2$. Therefore, our assumption of the existence of zeros of $W$ is wrong. Hence, we conclude that $ W$ does not vanish in $(a,b)$. This implies $W(x)<0$ for every $x \in(a,b)$, which in turn implies that $R$ is strictly decreasing on $(a,b)$.

Now, define a mapping $S : (\alpha, \beta) \to \mathbb{R}$ by $S(x) = {y_2(x)}/{y_1(x)},$ where $\alpha$ and $\beta$ are two consecutive zeros of $y_1$. Clearly, to show $S'(x)>0$ for every $x \in (\alpha, \beta)$, it is sufficient to show that $T(x)= y_2'(x)y_1(x)-y_1'(x)y_2(x) >0$ at each $x$ in $(\alpha, \beta)$. Notice that, again from \eqref{eq:y1coupledpart} and \eqref{eq:y2coupledpart} we have  $T(\alpha)>0$ and $T(\beta)>0$. Moreover, $T$ is a solution of a differential equation similar to \eqref{W-DE}. We claim that there are no zeros of $T$ in $(\alpha, \beta),$ and we can prove this by contradiction. Assume that for contradiction that there is $\kappa$ in $(\alpha, \beta)$ such that $T(\kappa)=0$. Then, with similar kinds of arguments as in the case of $W$, we show $T'(\kappa) \neq 0$. Then there must be finitely many zeros, and hence we must have the least zero, say $\kappa_1$, and its consecutive zero $\kappa_2$. Notice that in this case we have a sign of $T'(\kappa)$ opposite to the sign of $\kappa$, where $\kappa$ is a zero of $T$. Then, by continuity and the simplicity of the zeros, we have $T'(\kappa_1)<0$ and $T'(\kappa_2)>0$. Therefore $\kappa_1>0$ and $\kappa_2<0$, which is in contradiction to our construction of $\kappa_1<\kappa_2$. Therefore $T$ is also free of zeros. Hence we claim $T(x)>0$ for every $x \in (\alpha, \beta)$. Therefore, $S$ is strictly increasing. This completes the proof.
\end{proof}
	
\begin{lemma}\label{lemma:fDashNozero}
The derivative of the profile function $f'$ has no zeros in $(\tilde{R}, \tilde{R}_2].$ Moreover, if $k>0$ then $\tilde{R}_2<r_{1,N_1}$ and if $k<0$ then $\tilde{R}_2< r_{3,1}$.
\end{lemma}

\begin{proof}[\bf Proof]
We prove only the case when $k>0,$ since the argument for $k<0$ is quite similar. For the sake of contradiction, suppose that $f'(r_0) = 0$ for some $r_0 \in (\tilde{R}, \tilde{R}_2)$  and let $r_0$ be the first such zero. By using the relation
$$f'(r)=-\frac{\beta}{\sqrt{k}}\,\sin^{1-n/2}(\sqrt{k}r)\,\mathsf P^{-n/2}_{\gamma_{\beta,k,n}}(\cos(\sqrt{k}r))$$
for every $r \in (0, \pi/\sqrt{k})$ and the asymptotic formula \cite[eq. 14.8.1]{Nist}, we obsere that $f'(r)<0$ for $0<r<r_0$.
Now set $g(r) = -f'(r)$ and define $u = g\theta_1$ on the half-ball $D = \{(r, \theta) : 0 < r < r_0, \, \theta_1 > 0\}$. Here $(r, \theta)$ is the geodesic polar coordinate representation of a point in $M^{n}_{k},$ where $r$ is the geodesic distance and $\theta=(\theta_1,\dots,\theta_n)\in \mathbb{S}^{n-1}.$ Now, we observe that $u > 0$ in $D$, $u = 0$ on $\partial D$. For $k>0$, the Laplace-Beltrami operator is
$$\Delta_k = \partial_r^2 +(n-1)\sqrt{k}\cot(\sqrt{k}r)\,\partial_r+\frac{k}{\sin^2(\sqrt{k}r)}\,\Delta_{\mathbb S^{n-1}}.$$
From \cite[eq. 14]{SKAK24} we observe that $g$ satisfies
$$g''(r)+(n-1)\sqrt{k}\cot(\sqrt{k}r)g'(r)+\left(\beta-\frac{(n-1)k}{\sin^2(\sqrt{k}r)}\right)g(r)=0.$$ Note that $-\Delta_{\mathbb S^{n-1}}\theta_1=(n-1)\theta_1$. Hence we obtain
\begin{align*}
	\Delta_k u(r) &= \left[g''(r)+(n-1)\sqrt{k}\cot(\sqrt{k}r)g'(r)\right]\theta_1 +\frac{k}{\sin^2(\sqrt{k}r)}\,g(r)\,\Delta_{\mathbb S^{n-1}}\theta_1\\ &= \left[g''(r)+(n-1)\sqrt{k}\cot(\sqrt{k}r)g'(r) -\frac{(n-1)k}{\sin^2(\sqrt{k}r)}\,g(r)\right]\theta_1 =-\beta g(r)\theta_1.
\end{align*}
Thus, we have that $-\Delta_k u=\beta u$ in $D,$ $u=0$ on $\partial D,$ and $u>0$ in $D$. Hence $u$ is a positive eigenfunction on $D$, so $\lambda_1(D)=\beta$.
Since the hemisphere $H$ has exactly half the surface area of $\mathbb{S}^{n-1}$ we have that the volume $\operatorname{vol_k}(D)$ of $D$ is as follows
\begin{align*}
	\operatorname{vol_k}(D)&=\int_0^{r_0}\left(\frac{\sin(\sqrt{k}r)}{\sqrt{k}}\right)^{n-1}{\rm d}r \int_H {\rm d}\theta\\
	&=\frac12\,n\omega_n\int_0^{r_0}\left(\frac{\sin(\sqrt{k}r)}{\sqrt{k}}\right)^{n-1}{\rm d}r=\frac12 V(r_0).
\end{align*}
If $r_0<\tilde{R}_2$ then $\operatorname{vol_k}(D)=\frac12 V(r_0)<\frac12 V(\tilde{R}_2)=V(\tilde{R})$.
Now, we consider a geodesic ball $B(\rho)$ of radius $\rho$ with $V(\rho) = \operatorname{vol_k}(D)$. Then from the Faber-Krahn inequality \cite[p. 87]{Chavel} we obtain that
$\lambda_1(D) \ge \lambda_1(B(\rho)).$ Since $V(\rho) < V(\tilde{R})$ and $r \mapsto V(r)$ is strictly increasing, we have $\rho < \tilde{R}$, which means $B(\rho) \subseteq B(\tilde{R})$. Notice that $B(\tilde{R})-\overline{B(\rho)}$ is open. Therefore, in view of \cite[Corollary 1]{Chavel} we have $\lambda_1(B(\rho))>\lambda_1(B(\tilde{R}))$. This leads to the contradiction $$\beta=\lambda_1(D)\ge\lambda_1(B(\rho))>\lambda_1(B(\tilde{R}))=\beta.$$
It remains to exclude $f'(\tilde{R}_2)=0$. In this case, set $r_0=\tilde{R}_2$, and construct $D$ as above with this $r_0$. Then $\operatorname{vol_k}(D)=\frac 12 V(\tilde{R}_2)=V(\tilde{R}).$ Now, we again consider a ball $B(\rho)$ of radius $\rho$ such that $V(\rho)=\operatorname{vol_k}(D)$. From the Faber-Krahn inequality, we obtain that $\lambda_1(D) > \lambda_1(B(\rho))$. Here we have a strict inequality, as equality holds in the Faber–Krahn inequality if and only if $D$ is a geodesic ball, but $ D$ is not a geodesic ball. This also leads to a similar contradiction. Therefore $f'$ has no zeros in $(\tilde{R}, \tilde{R}_2]$. Finally, since $\tilde{R}<r_{1,N_1}$ we conclude that $\tilde{R}_2<r_{1,N_1}$.
\end{proof}

\begin{proof}[\bf Proof of Theorem \ref{thm:l-monotonicityk>0}]
By setting $x = \cos(\sqrt{k}r)$, we can rewrite $l$ in the $x$ domain as
$$l(r)=l\left(\arccos(x)/\sqrt{k}\right) = (1-x^2)^{\frac{n-1}{2}} [R(x)]^2$$
where $R(x)=y_1(x)/y_2(x)$ and  $x \in (-1,1)$. This implies
\begin{equation}\label{eq:l-monotonicityEQ1}
	l'(r) = \sqrt{k}(1-x^2)^{\frac{n-2}{2}}R(x) \left[ (n-1)x R(x) - 2(1-x^2)R'(x) \right]
\end{equation}
for every $r \in (0, \pi/\sqrt{k})$. Now, by using the difference-differential equations \eqref{eq:y1coupledpart} and \eqref{eq:y2coupledpart} along with the definition of $\gamma_{\beta,k,n}$ we observe that the derivative of the ratio $R'(x)$ satisfies the identity
$$(1-x^2)R'(x) = -\sqrt{1-x^2} + (n-1)x R(x) - \frac{\beta}{k}\sqrt{1-x^2}R^2(x)$$
for every $x \in (-1,1),$ which in turn implies that \eqref{eq:l-monotonicityEQ1} can be rewritten as
\begin{equation}
	l'(r) = \sqrt{k}(1-x^2)^{\frac{n-2}{2}}R(x) \left[ \sqrt{1-x^2} + \frac{\beta}{k}\sqrt{1-x^2}R^2(x) - (1-x^2)R'(x) \right]
\end{equation}
for every $r \in (0, \pi/\sqrt{k})$.

Observe that if $r \in (r_{2, N_2}, r_{1, N_1})$ then $x \in (x_{1,N_1}, x_{2, N_2})$. Consequently, in view of Theorem \ref{thm:FerrersRatio}, we have $R(x)<0$ and  $R'(x) < 0$ on this interval. Hence, the expression inside the bracket of \eqref{eq:l-monotonicityEQ1} is positive, and the expression outside is negative. Thus we obtain that $l'(r)<0$ for every $r \in (r_{2, N_2}, r_{1, N_1})$.

Now, we are going to show that $l$ is positive and strictly increasing in $(0, \tilde{R})$. Since $\mu=n/2 \geq 1$  we have $\lim_{x \to 1^{-}}R(x)=0$. From the definition of $W$, and equations \eqref{eq:y1coupledpart} and \eqref{eq:y2coupledpart} we obtain that
$$W(x_{2, N_2}) = -y_{2}'(x_{2, N_2})y_{1}(x_{2, N_2})= -\frac{\beta }{k}\frac{1}{\sqrt{1-x_{2, N_2}^2}}(y_1(x_{2, N_2}))^2<0.$$
From \cite[eq. 14.8.1 ]{Nist} we have that $W(x) \sim \frac{-1}{2n\left(\Gamma(n/2)\right)^{2}} \left(\frac{1-x}{2}\right)^{\frac{n-3}{2}}< 0$ as $x \to 1^{-}$. Now, through similar arguments as in the proof of Theorem \ref{thm:FerrersRatio}, we obtain that $W(x)<0$ for every $x \in (x_{2, N_2}, 1)$. Therefore $R'(x)<0$ for every $x \in (x_{2, N_2},1).$ Finally, from \cite[eq. 14.8.1 ]{Nist} we obtain that $R(x)>0$ in $(x_{2, N_2},1)$. Therefore $l'(r)>0$ on the corresponding interval $(0, r_{2,N_2})$.
\end{proof}
	
\begin{proof}[\bf Proof of Theorem \ref{thm:conicalZeros}]
Clearly from the series representations of conical functions (see \cite[p. 372]{Nist}), we have that if $\mu \neq 0,$ then $y_3(x)/y_4(x) \sim \frac{1}{\mu}\sqrt{\frac{x-1}{2}}$ as $x \searrow 1$. Therefore, $y_3(x)/y_4(x)$ is positive near $1$ whenever $\mu >0$ and negative whenever $\mu<0$. For $\mu=0$ we have $y_3(x) \to 1$ and $y_4(x) \sim \sqrt{\frac{x-1}{2}}$ as $x \to 1 ^{+}$. Therefore $y_3(x)/y_4(x)$ is also positive near $1$ whenever $\mu=0$.
	
Now, let us study the sign of $y_3(x)/y_4(x)$ at a zero of $ y_4 (x) $. Let $a$ and $c$ be points in $(1,\infty)$ such that $y_4(a)=0$ and $y_3(c)=0$.
We observe that $y_3$ and $y_4$ have the Taylor series expansions
$$ y_3(x)= y_3(a)+y_3'(a)(x-a)+\frac{y_3''(a)}{2} (x-a)^2+\ldots,$$
$$ y_4(x)= y_4'(a)(x-a)+ \frac{y_4''(a)}{2}(x-a)^2+\ldots$$
for every $x \in (a-\delta, a+\delta)$ for some $\delta>0$. Thus $\frac{y_3(x)}{y_4(x)} \sim \frac{y_3(a)}{y_4'(a)(x-a)}$ as $x \to a$. Similarly, by using the Taylor series around $c$ we obtain $\frac{y_3(x)}{y_4(x)} \sim \frac{y_3'(c)}{y_4(c)}(x-c)$ as $x \to c$. By using \eqref{eq:y4coupledpart} we can write
\begin{equation}
\frac{y_3(x)}{y_4(x)} \sim \frac{-\sqrt{a^2-1}}{((\mu-1/2)^2+\tau^2)(x-a)} \quad \text{as } x \to a
\end{equation}
and similarly from \eqref{eq:y3coupledpart} we can write
$$ \frac{y_3(x)}{y_4(x)} \sim \frac{x-c}{\sqrt{c^2-1}} \quad \text{as } x \to c.$$
Therefore $y_3(x_{4,1}^-)/y_4(x_{4,1}^-)>0$, $y_3(x_{3,1}^{-})/y_4(x_{3,1}^{-})<0$ and $y_3(x_{3,1}^{+})/y_4(x_{3,1}^{+})>0$.
Now, we claim that for $\mu\geq 0$ we have $x_{4,1}<x_{3,1}$. For the sake of contradiction suppose that $x_{3,1}\leq x_{4,1}.$ In this case $y_3(x)/y_4(x)$ is positive on $(1, x_{3,1})$ and changes sign from positive to negative at $x_{3,1}$. Since $y_3(x_{4,1}^-)/y_4(x_{4,1}^-)>0$ then there must be another zero of $y_3$ before $x_{4,1}$. This is in contradiction with the interlacing properties of zeros (see Lemma \ref{lemma:interlace}). Similarly, for $\mu<0$ we have $x_{3,1}<x_{4,1}$.
\end{proof}

\end{document}